\documentstyle[11pt]{article}
\normalbaselines

\newtheorem{Theorem}{Theorem}[section]
\newtheorem{Definition}[Theorem]{Definition}
\newtheorem{Proposition}[Theorem]{Proposition}
\newtheorem{Lemma}[Theorem]{Lemma}
\newtheorem{Corollary}[Theorem]{Corollary}
\newtheorem{Remark}[Theorem]{Remark}

\newcounter{alph}

\title{On the fundamental group of an abelian cover
\thanks{Both authors were supported by the
SCIENCE Program (Contract n. SCI-0398-C(A)).}\thanks{A.M.S.
classification 14E20}}
\author{Rita Pardini\thanks{A member of
G.N.S.A.G.A. of C.N.R.}, Francesca Tovena}
\date{}
\begin{document}
\maketitle
\def\qed{\hbox{\hskip 6pt\vrule width6pt height7pt depth1pt \hskip
1pt}\bigskip}

{\bf R\'esum\'e.} Soient $X$ et $Y$ deux vari\'et\'es complexes
projectives et lisses de dimension $n\geq 2$ et soit ${f: Y\to X}$
un rev\^etement abelien totalement ramifi\'e.
Alors l'application ${f_*:\pi_1(Y)\to \pi_1(X)}$ est
surjective et donne une extension centrale:
\begin{equation}\label{***}
0\to K\to \pi_1(Y)\to\pi_1(X)\to 1
\end{equation}
o\`u $K$ est un groupe fini.

Nous montrons comment le noyau $K$ et la classe de cohomologie
$c(f)\,\in\, H^2(\pi_1(X),K)$ de (\ref{***}) peuvent \^etre
calcul\'es en termes de classes de Chern des composantes du
diviseur critique de $f$ et des sous-faisceaux inversibles de
$f_*{\cal O}_Y$ stables sous l'action du groupe de Galois.

{\bf Abstract.} Let $X$, $Y$ be smooth complex projective varieties of
dimension $n\geq 2$ and let ${f: Y\to X}$  be a totally ramified
abelian cover. Then the map ${f_*:\pi_1(Y)\to \pi_1(X)}$ is
surjective and gives rise to a central extension:
\begin{equation}\label{**}
0\to K\to \pi_1(Y)\to\pi_1(X)\to 1
\end{equation}
where $K$ is a finite group.

Here we show how the kernel $K$ and the cohomology class $c(f)\,\in\,
H^2(\pi_1(X),K)$ of (\ref{**}) can be computed in terms of the Chern
classes of the components of the branch divisor of $f$ and of the
eigensheaves of $f_*{\cal O}_Y$ under the action of the Galois group.

\section{Introduction.}
\hspace{6 mm}This work generalizes a result of Catanese and the second
author, who analyze in \cite{kn:Cato} the fundamental group of a
special type of covering ${f:Y\rightarrow
X}$, with Galois group $({\bf Z}/m{\bf
Z})^{2}$, of a complex smooth projective surface $X$, the so-called
"$m$-th root extraction" of a divisor $D$ on $X$.

By means of standard topological methods, the fundamental group
$\pi_{1}(Y)$ can be described in that case
as a central extension of the group $\pi_{1}(X)$, as follows:
\begin{equation}
\label{ext1}
0\longrightarrow {\bf Z}/r{\bf Z}\longrightarrow \pi_{1}(Y)
\longrightarrow \pi_{1}(X) \longrightarrow 1,
\end{equation}
$r$ being a divisor of $m$ which depends only on the divisibility of
$\pi^{*}(D)$ in $H^2(\tilde{X},{\bf Z})$, where
${\pi:\tilde{X}\rightarrow X}$ is the universal covering of $X$.

The main result of \cite{kn:Cato} (see Thm.2.16) is that
the group cohomology class corresponding to the extension
(\ref{ext1}) can be explicitly computed in terms of the first
Chern class of $D$.

This is an instance of a more general philosophy: in principle,
it should be possible to recover all the information
about an abelian cover ${f:Y\to X}$ from the "building data"
of the cover, i.e., from the Galois group $G$, the
components of the branch locus, the inertia subgroups
and the eigensheaves of ${f_{*}{\cal O}_{Y}}$ under the natural
action of $G$ (see section 2 or \cite{kn:Rita} for more
details).

Actually, the description of the general abelian cover given in
\cite{kn:Rita} enables us  to treat (under some mild assumptions on the
components $D_{1},\ldots D_{k}$ of the branch locus) the case of any
totally ramified abelian covering ${f:Y\rightarrow X}$, with $X$ a
complex projective variety of dimension at least $2$ (cf. section 2
for the definition of a totally ramified abelian cover).

Using the same methods as in \cite{kn:Cato}, we show that
$\pi_{1}(Y)$  is a central extension as before:
\begin{equation}
\label{ext2}
0\longrightarrow K\longrightarrow \pi_{1}(Y) \longrightarrow
\pi_{1}(X) \longrightarrow 1\,,
\end{equation}
where $K$ is a finite abelian group which is determined by the
building data of the cover and the cohomology classes of
$\pi^{*}(D_{1}),\ldots \pi^{*}(D_{k})$ (cf. Prop.\ref{top}).

Consistently with the above "philosophy", a statement analogous to
Thm.2.16 of \cite{kn:Cato} actually holds in the general case: our
main result (Thm.\ref{mtgen}, Rem.\ref{dipdaLchi}) can be summarized by
saying that the class of the extension (\ref{ext2}) can be recovered
from the Chern classes of the $D_j$'s and of the eigensheaves of
$f_*{\cal O}_Y$; in some special case, this relation can be described
in a particularly simple way (Thm.\ref{mt}, Cor.\ref{mtcor},
Rem.\ref{rem}). Moreover, one can construct examples of not
homeomorphic varieties realized as covers of a projective variety $X$ with the
same Galois group, branch locus and inertia subgroups (cf.
Rem.\ref{dipdaLchi}).

The idea of the proof is to exploit a natural representation of
$\pi_{1}(Y)$ on a vector bundle on the universal covering
$\tilde{X}$ of $X$ and the spectral sequence describing the
cohomology of a quotient, in order to relate the group
cohomology class of the extension (\ref{ext2}) to
the geometry of the covering. These are basically
the same ingredients as in the proof of \cite{kn:Cato},
but we think that we have reached here a more conceptual
and clearer understanding of the argument.

\noindent
{\bf Acknowledgements:} we wish to express our heartfelt
thanks to Fabrizio Catanese, who suggested that the result
of \cite{kn:Cato} was susceptible of generalization
and encouraged us to investigate this problem.

\section{A brief review of abelian covers.}

\hspace{6 mm}In this section we set the notation and, for the reader's
convenience, we collect here the definitions and the
notions concerning abelian covers that will be needed later. For
further details and proofs, we refer to \cite{kn:Rita},
sections 1 and 2.

Let $X$, $Y$ be complex algebraic varieties of dimension at
least 2, smooth and projective, and let ${f:Y\to X}$
be a finite abelian cover, i.e. a Galois cover with finite
abelian Galois group $G$.

The bundle $f_{*}({\cal O} _{Y})$ splits as a sum of one dimensional
eigensheaves under the action of $G$, so that one has:
\begin {equation}
\label{splitting}
f_{*}({\cal O} _{Y})=
\bigoplus_{\chi \in G^{*}}L_{\chi}^{-1}
={\cal O} _{X}
\oplus ( \bigoplus _{\chi \in G^{*}\setminus \{1\}}
L_{\chi}^{-1})
\end{equation}
where $G^{*}$ denotes the group of characters of G and $G$ acts on
$L_{\chi}^{-1}$ via the character $\chi$.

We warn the reader that the notation here and in the next section
is dual to the one adopted in \cite{kn:Rita}; however this
does not affect the formulas quoted from there.

Under our assumptions, the ramification locus of $f$ is a divisor.
Let $D_{1}, \ldots D_{k}$ be the irreducible components of
the branch locus $D$ and let $R_{j} = f^{-1}(D_{j})$,
${j = 1,} \ldots k$. For every ${j = 1, \ldots k}$, one defines the
{\em inertia subgroup} $G_{j} = \{g\in G | g(y)=y \;
{\rm for \;each}\;y\in R_{j}\}$. Given any point
$y_{0} \in R_{j}$, one obtains a natural representation
of $G_{j}$ on the normal space to $R_{j}$ at $y_{0}$
by taking differentials. The corresponding character, that we
 denote by $\psi_{j}$, is independent of the choice
of the point $y_{0} \in R_{j}$. By standard results, the subgroup
$G_{j}$ is cyclic and the character $\psi_{j}$ generates the group
$G_{j}^{*}$ of the characters of $G_{j}$. We denote by $m_{j}$ the
order of $G_{j}$, by $m$ the least common multiple of the $m_{j}$'s
and by $g_{j}$ the generator of $G_{j}$ such that
$\psi_{j}(g_{j}) = exp(\frac{2\pi \sqrt{-1}}{m_{j}})$.

In what follows we will always assume that the cover
${f:Y\rightarrow X}$ is {\em totally ramified}, i.e., that
the subgroups $G_{j}$ generate $G$; then the group of
characters $G^{*}$ injects in
${\bigoplus_{j=1}^{k}G_{j}^{*}}$ and every
${\chi \in G^{*}}$ may be written uniquely as:
\begin{equation}
\label{char0}
\chi = \sum_{j=1}^{k} \;a_{\chi,j}\,\psi_{j},
\hspace{10mm}\,0\leq a_{\chi,j}< m_{j} \;{\rm for\,every}\,j\,.
\end{equation}
In particular, let $\chi_{1}, \ldots \chi_{n} \in G^{*}$
be such that $G^{*}$ is the direct sum of the cyclic
subgroups generated by the $\chi $'s, and let $d_{i}$
be the order of $\chi_{i}$, $i = 1, \ldots n$. Write:
\begin{equation}
\label{char}
\chi_{i} = \sum_{j=1}^{k} \;a_{ij}\,\psi_{j},
\;\;\,\;0\leq a_{ij}< m_{j}\,, \;\;\, i= 1, \ldots n\,.
\end {equation}
Then one has (\cite{kn:Rita}, Prop.2.1):
\begin{equation}
\label{eqstr}
d_{i} L_{\chi_{i}}
\equiv \sum_{j=1}^{k}\frac{d_{i}a_{ij}}{m_{j}}\;D_{j}
\;\;\;\;i= 1, \ldots n\,
\end{equation}
the corresponding isomorphism of line bundles being induced
by multiplication in the ${\cal O}_{X}$-algebra
$f_{*}{\cal O}_{Y}$. More generally, if $\chi=\sum_{i=1}^n
b_{\chi,i}\chi_i$, with $0\leq b_{\chi,i}<d_i
\;\forall\,i$, one has (\cite{kn:Rita}, Prop.2.1):
\begin{equation}\label{Lchigen}
L_{\chi}\equiv \sum_{i=1}^n b_{\chi,i}L_{\chi_{i}} -
\sum_{j=1}^k q^{\chi}_{j} D_j\,.
\end{equation}
where $q^{\chi}_{j}$ is the integral part of the rational
number ${\sum_{i=1}^n \frac{b_{\chi,i}a_{ij}}{m_j}}$, $j=1,\ldots
k$.

Equations (\ref{eqstr}) are the characteristic relations of an
abelian cover. Actually, since $X$ is complete, for assigned $G$,
$D_{j}$, $G_{j}$, $\psi_{j}$,  $j = 1, \ldots k$, to each set of line
bundles $L_{\chi_{i}}$, $i = 1, \ldots n$, satisfying (\ref{eqstr})
there corresponds a unique, up to isomorphism, $G$-cover of $X$,
branched on the $D_{j}$'s and such that $G_{j}$ is the inertia
subgroup of $D_{j}$ and $\psi_{j}$ is the corresponding character
(\cite{kn:Rita}, Thm.2.1). Moreover, the cover is actually smooth
under suitable assumptions on the building data.

\section{The fundamental group and the universal covering of $Y$.}

\hspace{6 mm}We keep the notation introduced in the previous section.
\begin{Definition}{\rm (\cite{kn:MM}, pag.218)}
A smooth divisor $\Delta$ on a variety $X$ is called
{\rm flexible}
if there exists a smooth divisor $\Delta ' \equiv \Delta$
such that $\Delta ' \cap \Delta \neq \emptyset$
and $\Delta$ and $\Delta '$ meet transversely.
\end{Definition}
We recall that a flexible divisor on a projective surface is connected
(see \cite{kn:Ca}, Remark 1.5). Hence, by considering a general linear
section, one deduces easily that a flexible divisor on a
projective variety of dimension $\geq 2$ is connected.
\begin{Proposition}
\label{top}
Let $X$, $Y$ be smooth projective varieties over {\bf C} of
dimension $n$ at least 2. Let ${f:Y\to X}$ be a totally ramified
abelian cover branched on irreducible, flexible and ample divisors
$\{D_j\}_{j=1,\ldots k}$. Then:
\begin{list}%
{\alph{alph})}{\usecounter{alph}}
\item
The natural map ${f_{*}:\pi_{1}(Y)\to \pi_{1}(X)}$
is surjective.
\item Let $K = ker(f_{*})$; then $K$ is finite and
\begin{equation}
\label{ext}
0\longrightarrow K\longrightarrow \pi_{1}(Y) \longrightarrow
\pi_{1}(X) \longrightarrow 1.
\end{equation}
is a central group extension.
\item Let ${\pi :\tilde{X}\to X}$ be the universal covering of $X$
and ${\tilde{D} = \pi^{-1}(D)}$;
then $\tilde{D}_{j}=\pi^{-1}(D_{j})$ is connected, j = 1, \ldots k.
Denote by $H^{i}_{c}$ the cohomology with compact supports and by
${\rho : H^{2n-2}_{c}(\tilde{X}) \to H^{2n-2}_{c}(\tilde{D})
\cong \bigoplus ^{k}_{j=1}\,{\bf Z}\tilde{D}_{j}}$ the restriction
map. Finally, let $\sigma$ be the map defined by:
\begin{equation}
\begin{array}{crclc}
\sigma:&H^{2n-2}_{c}(\tilde{D})
\cong \bigoplus ^{k}_{j=1}\,{\bf Z}\tilde{D}_{j}
&\to& \bigoplus_{j=1}^{k}G_{j}&\\
&{\tilde{D}_{j}} &\mapsto &g_{j}&.
\end{array}
\end{equation}

Then ${N = \ker(\bigoplus G_{j} \to G)}$ contains
${\rm Im}(\sigma \circ \rho)$ and $K$ is isomorphic
to the quotient group $N/{\rm Im}(\sigma \circ \rho)$.
\end{list}
\end{Proposition}

{\sc Proof.} a) and the fact that the extension (\ref{ext}) is central
can be proven exactly as in  \cite{kn:Ca}, Thm.1.6 and in
\cite{kn:Cato}, Lemma 2.1.

For the proof of c) (that implies that $K$ is finite), we refer the
reader to \cite{kn:Ca}, Prop.1.8 and to \cite{kn:Cato}, proof of
Thm.2.16, Step I. One only has to notice that, by Lefschetz theorem
(cf. \cite{kn:Bott}, Cor. of Thm.1), $\pi_{1}(D_{j})$ surjects onto
$\pi_{1}(X)$, hence $\tilde{D}_{j}=\pi^{-1}(D_{j})$ is connected and
smooth for every $j = 1, \ldots k$.\hfill\qed

\begin{Remark}{\rm
\begin{list}%
{\alph{alph})}{\usecounter{alph}}
\item
{}From Prop.\ref{top}, c), it follows in particular
that the kernel $K$ of the surjection ${f_*:\pi_1(Y)\to \pi_1(X)}$
does not depend on the choice of the solution $L_\chi$ of
(\ref{eqstr}), once $G$, the $g_j$'s and the class of the
$\tilde{D}_j$'s in $H^2(\tilde{X},{\bf Z}/m_j{\bf Z})$, $j=1,\ldots
k$, are fixed.
\item
If $f:Y\to X$ is an abelian cover
as in the hypotheses of Prop.\ref{top}, then
$H^1(Y,{\cal O}_Y)\cong H_1(X,{\cal O}_X)$ by (\ref{splitting}) and the
Kodaira Vanishing Theorem. Moreover, according to Prop.\ref{top}, a)
the map ${f_*:H_1(Y,{\bf Z})\to H_1(X,{\bf Z})}$ is surjective; thus
the map ${f_*: alb(Y)\to alb(X)}$ between the Albanese varieties is an
isomorphism.
\end{list}
 }
\end{Remark}

\begin{Proposition}
\label{tilde}
In the same hypotheses as in Prop.\ref{top}, let ${q:\tilde{Y}\to Y}$
be the universal cover of $Y$ and let
${\tilde{f\,}:\tilde{Y}\to \tilde{X}}$ be the map
lifting $f:Y\to X$. Then $\tilde{f\,}$ is a totally ramified
abelian cover of $\tilde{X}$ with group
${\tilde{G} = (\bigoplus_{j=1}^k G_{j})/{\rm Im}(\sigma \circ \rho)}$,
branched on $\tilde{D}$.
\end{Proposition}

{\sc Proof.} By diagram chasing, it is easy to show that
${\pi_{1}(\tilde{X}\setminus \tilde{D})}$ is isomorphic
to the kernel $V$ of the surjection
${\pi_{1}(X\setminus D)\to \pi_{1}(X)}$ induced by the inclusion
${X\setminus D\subset X}$.
Since the $D_{i}$'s are flexible, one proves as in (\cite{kn:Cato},
Lemma 2.1) that $V$ is an abelian group. It follows that
$\tilde{f\,}$, being branched on $\tilde{D}$, is an abelian cover.

Consider now the fiber product $Y'$ of ${f:Y\to X}$ and
${\pi:\tilde{X} \to X}$, with the natural maps
${f':Y'\to\tilde{X}}$  and ${q':Y'\to Y}$; $f'$ is a $G$-cover
ramified on $\tilde{D}$ and $q'$ is unramified. According to
Prop.\ref{top}, b), the universal covering ${q:\tilde{Y}\to Y}$ of
$Y$ factors as ${q=q' \circ q''}$, for a suitable unramified cover
${q'':\tilde{Y} \to Y'}$ with group $K$, giving a commutative diagram
as follows:   \begin{equation} \label{diagramma}
\begin{picture}(60,85)(-4,-40)
\put (8,0){\vector(1,0){20}}
\put (9,30){\vector(1,-1){20}}
\put (0,30){\vector(0,-1){19}}
\put (0,-8){\vector(0,-1){19}}
\put (42,-8){\vector(0,-1){19}}
\put (8,-38){\vector(1,0){20}}
\put (-4,-2){$Y'$}
\put (36,-2){$\tilde{X}$}
\put (-3,36){$\tilde{Y}$}
\put (-4,-40){$Y$}
\put (36,-40){$X$}
\put(23,20){${\scriptstyle \tilde{f\,}}$}
\put(13,3){${\scriptstyle f'}$}
\put(-10,20){${\scriptstyle q''}$}
\put(13,-35){${\scriptstyle f}$}
\put(-10,-20){${\scriptstyle q'}$}
\put(44,-20){${\scriptstyle \pi}$}
\end{picture}
\end{equation}
In particular, ${K\cong \pi_{1}(Y')}$ and
${\tilde{f\,}=f'\circ q''}$. Hence, the Galois group $\tilde{G}$ of
$\tilde{f\,}$ is given as an extension:
\begin{equation}
\label{extgrouO}
0\longrightarrow K \longrightarrow\tilde{G}
\longrightarrow G \longrightarrow 0
\end{equation}
Moreover, if one denotes by $\tilde{G}_{j}$ the inertia subgroup of
$\tilde{D}_{j}$
with respect to $\tilde{f\,}$, then $\tilde{G}_{j}$ maps
isomorphically onto $G_{j}$ for every $j=1, \ldots k$.
The isomorphism
${\tilde{G} = (\bigoplus G_{j})/{\rm Im}(\sigma \circ \rho)}$
can be obtained by computing the fundamental group of
$\tilde{Y}$ as in \cite{kn:Cato}, proof of Thm.2.16.
\hfill\qed

\vspace{3mm}
The following lemma will be used in the next section.
\begin{Lemma}
\label{commut}
Consider the subgroups $\pi_{1}(Y)$ and $\tilde{G}$ of
${Aut(\tilde{Y})}$; then one has:
\begin{equation}
\beta g = g\beta\,\hspace{15mm}\forall g\in \tilde{G},
\forall\beta \in \pi_{1}(Y).
\end{equation}
\end{Lemma}

{\sc Proof.} Since the cover $\tilde{f\,}$ is totally ramified, it is
enough to show that all the elements of $\pi_{1}(Y)$ commute with
$g_{j}$,  $j=1, \ldots k$.

Let ${\beta \in \pi_{1}(Y)}$ and fix $j=1, \ldots k$.
We remark firstly that ${\beta g_{j}\beta^{-1}}$ is actually an
element of ${\tilde{G}\subset Aut(\tilde{Y})}$. In fact, consider
the classes represented by ${\beta g_{j}}$ and ${g_{j}\beta}$ modulo
$K$: they do coincide as automorphisms of ${Y'\subseteq
Y\times\tilde{X}}$, since the group ${G\times\pi_{1}(X)}$ acts there
via the natural action on the components. So, ${\beta
g_{j}\beta^{-1}g_{j}^{-1}\in K}$ and  ${\beta g_{j}\beta^{-1}\in
g_{j}K\subseteq\tilde{G}}$, as desired.

By diagram (\ref{diagramma}), we have
$\tilde{R}_{j}=\tilde{f\,}^{-1}(\tilde{D}_{j})=q^{-1}(R_{j})\;\forall$
$j=1, \ldots k$. Since $R_{j}=f^{-1}(D_{j})$ is ample and connected,
the same argument as in the proof of Lemma \ref{top}, c) shows
that $\tilde{R}_{j}$ is connected. So,
${\beta \tilde{R}_{j}=\tilde{R}_{j}}$ and ${\beta g_{j}\beta^{-1}}$
fixes ${\tilde{R}_{j}}$ pointwise, namely
${\beta g_{j}\beta^{-1}\in \tilde{G}_{j}}$.

Finally, recalling the definition of the character ${\psi_{j}\in
G_{j}^{*}}$ introduced in section 2, one checks immediately that
${\psi_{j}(\beta
g_{j}\beta^{-1})=\psi_{j}(g_{j})}$.
The conclusion now follows from the faithfulness of $\psi_{j}$.

\hfill\qed

\section{Computing the cohomology class of the central extension
$0\to K \to\pi_{1}(Y)\to \pi_{1}(X)\to 1$.}

\hspace{6 mm}We keep the notation and the assumptions introduced
in the previous sections, unless the contrary is explicitly
stated. We need two technical Lemmas in order to state the main result
of this paper.

\begin{Lemma}
\label{Lemma}
Let ${\cal H}$ be a finite abelian group and $\zeta_1$,\ldots
$\zeta_m\,\in\,{\cal H}$ be such that ${\cal H}=\bigoplus_{j=1}^m
<\zeta_j>$ is the direct sum of the cyclic
subgroups generated by $\zeta_j$, j=1, \ldots m; denote
by $h_j$ the order of $\zeta_j$. Let ${p\,\in\,{\bf Z}}$ be a prime
and ${\cal H}_{p}$ be the p-torsion subgroup of ${\cal H}$. Let
$\chi_{1}, \ldots \chi_{t} \in {\cal H}_{p}$ such that
${<\chi_{1}, \ldots \chi_{t}>=\bigoplus_{i=1}^t<\chi_{i}>}$. Finally,
let $d_{i}$ be the order of $\chi_{i}$ and write
${\chi_{i}=\sum_{j=1}^{h} a_{ij}\zeta_{j}}$ with ${0\leq a_{ij}<h_j}$.

Then, $\forall$ $x_{1}, \ldots x_{t} \in {\bf Z}$ and $\forall$
$\gamma \geq 1$, the system:
\begin{equation}
\label{sistlemma} \sum_{j=1}^{m}\frac{d_{i}a_{ij}}{h_{j}}s_{j}\equiv
x_{i}\;\;\; {\rm mod}\,p^{\gamma}\;\;\;\;\;\;\;\;\;i=1,\ldots t
\end{equation}
admits a solution $(s_{1}, \ldots s_{m})\,\in\,{\bf Z}^m$.
\end{Lemma}

{\sc Proof.} We set $c_{ij}=\frac{d_{i}a_{ij}}{h_{j}}$ and, for
$x\,\in\,{\bf Z}$, we denote by $\overline{x}$ the class of $x$ in
${\bf Z}/p{\bf Z}$. We proceed by induction on $\gamma$.

Let $\gamma =1$. We show that
the matrix $(\overline{c}_{ij})$ has rank $t$.

Let $y_{1}, \ldots y_{m} \in {\bf Z}$ and assume that:
\begin{equation}
\sum_{i}\overline{c}_{ij}\overline{y}_{i}=0\;\;\;\;\;\forall\,j=1,
\ldots m\,.
\end{equation}
This implies that:
\begin{equation}
\sum_{i=1}^tc_{ij}y_{i}\equiv 0\;\;\;{\rm
mod}\,p\;\;\;\;\;\;\;\forall\,j=1, \ldots m
\end{equation}
so that:
\begin{equation}
\sum_{i=1}^t\frac{y_{i}d_{i}}{p}\frac{a_{ij}}{h_{j}}\;\in\,
{\bf Z}\;\;\;\;\forall\,j=1, \ldots m\,.
 \end{equation}
Recalling that $p$ divides $d_{i}$ $\forall\,i$, we deduce that:
\begin{equation}
\sum_{i=1}^t\left(\frac{y_{i}d_{i}}{p}\right)\,a_{ij}\equiv 0
\;\;\;{\rm mod}\,h_{j}\;\;\;\;\;\;\;\forall\,j=1, \ldots m\,,
\end{equation}
so that $\sum_{i}\frac{y_{i}d_{i}}{p}\chi_{i}$ is the zero element in
${\cal H}$. By the hypothesis on the $\chi_{i}$'s, it
follows that:
\begin{equation}
\frac{y_{i}d_{i}}{p}\equiv 0\;\;\;{\rm mod}\,d_{i}
\end{equation}
and finally:
\begin{equation}
y_{i}\equiv 0\;\;\;{\rm mod}\,p
\end{equation}
showing, as desired, that the rows of the matrix $(\overline{c}_{ij})$
are linearly independent over ${\bf Z}/p{\bf Z}$.

Let now $\gamma > 1$ and assume by inductive hypothesis that ${(s_{1},
\ldots s_{m})\in {\bf Z}^{m}}$ is a solution of the system
(\ref{sistlemma}).

We set $s_{j}'=s_{j}+\delta_{j}p^{\gamma}$ and we look for a suitable
choice of the integers $\delta_{j}$. We have:
\begin{equation}
\begin{array}{ccl}
\sum_{j}c_{ij}s_{j}'&=&\sum_{j}c_{ij}s_{j}+p^{\gamma}
\sum_{j}c_{ij}\delta_{j}\\
&=&x_{i}+p^{\gamma}y_{i}+p^{\gamma}
\sum_{j}c_{ij}\delta_{j}\;\;\;\;\;\exists\,y_{i}\,\in\,{\bf Z},
\end{array}
\end{equation}
so that:
\begin{equation}
\sum_{j}c_{ij}s_{j}'\equiv x_{i}\;\;{\rm mod}\,p^{\gamma +1}
\;\;\;\Longleftrightarrow\;\;\;\sum_{j}c_{ij}\delta_{j} \equiv -y_{i}
\;\;\;{\rm mod}\,p
\end{equation}
and the latter system has a solution, by the case $\gamma = 1$. This
conclude the proof.\hfill\qed

We come back to the study of the cover $f$:
\begin{Lemma}
\label{coeffinK}
Let $A$ be the subgroup of $Pic(X)$ generated by
$D_{1}, \ldots D_{k}$ and $L_{\chi}$, $\chi\in G^{*}$. Then there
exist $M_{1}, \ldots M_{q}\in Pic(X)$ such that
${A=\bigoplus_{l=1}^{q}<M_{l}>}$ and
\begin{equation}
\left(\begin{array}{c}
D_{1}\\
\vdots\\
D_{k}
\end{array}\right)\;\;\equiv\;\;C\,
\left(\begin{array}{c} M_{1}\\
\vdots\\
M_{q}
\end{array}\right)
\end{equation}
 where $C=(c_{jl})$ is a matrix with integral coefficients such
that each column $(c_{jl})_{j=1,\ldots k}$ represents an element of
${N=\,{\rm ker}\,(\bigoplus_{j=1}^kG_j\to G)}$.
\end{Lemma}

{\sc Proof.} $A$ is a finitely generated abelian group, so one can
write ${A=F\bigoplus T}$, where $T$ is the torsion part of $A$ and $F$
is free.

Denote by $\{\xi_{l}\}_{l}$ a set of free generators of $F$ and by
$\{\eta_{l}\}_{l}$ a set of generators of $T$ such that
${T=\bigoplus<\eta_{l}>}$ and the order $o(\eta_{l})$ of $\eta_{l}$ is
the power of a prime, $\forall l$. Let finally $\chi_{i}$ be
generators of $G^{*}$ such that $G^{*}=\bigoplus_{i=1}^n<\chi_{i}>$
and the order $o(\chi_{i})$ of $\chi_{i}$ is the power of a prime,
$\forall i$.

One can write:
\begin{eqnarray}
\label{Lchii}L_{\chi_{i}}\equiv
\sum_l\lambda_{il}\eta_{l}+\sum_l\lambda_{il}'\xi_{l}\hspace{10mm}
\forall\,i=1,\ldots n,\\
\label{Djlemma}D_{j}\equiv\sum_l c_{jl}\eta_{l}+\sum_l
c_{jl}'\xi_{l}\hspace{10mm} \forall\,j=1,\ldots k,
\end{eqnarray}
where the coefficients $\lambda_{il}'$ and $c_{jl}'$ are
uniquely determined, whereas $\lambda_{il}$ and $c_{jl}$ are
determined only up to a multiple of $o(\eta_{l})$.

We can apply the analysis of section 2 to the cover $f$. We write
${\chi_i=\sum_{j=1}^{k}a_{ij}\psi_j}$, with ${0\leq
a_{ij}<m_j}$, and we set ${d_{i}=o(\chi_i)}$ as in the previous
Lemma; the equations (\ref{eqstr}) become here:
\begin{equation} \label{eqstrf}
d_{i} L_{\chi_{i}}
\equiv \sum_{j=1}^{k}\frac{d_{i}a_{ij}}{m_{j}}\;D_{j}\;
\hspace{10mm}i=1,\ldots n,
\end{equation}
so that we must have:
\begin{equation}
d_{i}\lambda_{il}'=\sum_{j=1}^k\frac{d_{i}a_{ij}}{m_{j}}\;c_{jl}'
\hspace{10mm}i=1,\ldots n,
\end{equation}
showing that $(c_{jl}')_{j=1,\ldots k}$ represents an element of $N$,
$\forall \,l$: in fact, by duality, ${(t_1,\ldots,t_k)\,\in\,
{\bf Z}^k}$ represents an element of $N$ if and only if it satisfies
the relations:
\begin{equation}\label{relN}
\sum_{j=1}^k\frac{a_{ij}}{m_{j}}t_j\;\in\;{\bf
Z}\hspace{15mm}\forall\,i=1,\ldots n. \end{equation}

For the coefficients of the torsion part, we have:
\begin{equation}
d_{i}\lambda_{il}\,\eta_{l}=\left(
\sum_{j=1}^k\frac{d_{i}a_{ij}}{m_{j}}\;c_{jl}\right)\;\eta_{l}
\end{equation}
so that:
\begin{equation}
\label{cong}
d_{i}\lambda_{il}\equiv
\sum_{j=1}^k\frac{d_{i}a_{ij}}{m_{j}}\;c_{jl}
\;\;\;{\rm mod}\,o(\eta_{l})\;.
\end{equation}

We fix an index $l$. Let $p$ be a prime such that
$o(\eta_l)=p^{\alpha}$. We want to show that, for a suitable
choice of the $c_{jl}$, the following relation holds $\forall\,i=1,
\ldots n$:
\begin{equation}
\label{tesi}
d_{i}\lambda_{il}\equiv
\sum_{j=1}^k\frac{d_{i}a_{ij}}{m_{j}}\;c_{jl}
\;\;\;{\rm mod}\,d_{i}\,.
\end{equation}
Let $\chi_{i}$ be a generator such
that ${d_{i}\equiv 0}$ mod $p$ and set
${d_{i}=p^{\alpha_{i}}}$.
By (\ref{cong}), it is enough to consider the case in which
${\alpha<\alpha_{i}}$.

Setting $c_{jl}''= c_{jl}+p^{\alpha}s_{j}$ and recalling (\ref{cong}),
one has:
\begin{equation}
\begin{array}{ccl}
\sum_{j=1}^k\frac{d_{i}a_{ij}}{m_{j}}\;c_{jl}''&=&
\sum_{j=1}^k\frac{d_{i}a_{ij}}{m_{j}}\;c_{jl}+
p^{\alpha}\sum_{j=1}^k\frac{d_{i}a_{ij}}{m_{j}}\;s_{j}\\
&=&d_{i}\lambda_{il}-p^{\alpha}x_{i}+p^{\alpha}
\sum_{j=1}^k\frac{d_{i}a_{ij}}{m_{j}}\;s_{j}
\end{array}
\end{equation}
for a suitable choice of integers $x_{i}$. One concludes that the
relation (\ref{tesi}) holds if and only if:
\begin{equation}
\sum_{j=1}^k\frac{d_{i}a_{ij}}{m_{j}}\;s_{j}\equiv x_{i}
\;\;\;{\rm mod}\,p^{\alpha_{i}-\alpha}\;.
\end{equation}

Let $\beta = {\rm max}\,\{\alpha_{i}-\alpha\}_{i}$. The system of
congruences:
\begin{equation}
\sum_{j=1}^k\frac{d_{i}a_{ij}}{m_{j}}\;s_{j}\equiv x_{i}
\;\;\;{\rm mod}\,p^{\beta}\;\hspace{10mm}\forall\,i \mbox{ such
that } d_{i}\equiv 0\;\mbox{ mod } p
\end{equation}
admits a solution by Lemma \ref{Lemma}. So, we can assume that the
coefficients $(c_{jl})_{j=1,\ldots k}$ in (\ref{Djlemma}) satisfy
(\ref{tesi}) for every $i$ such that $d_{i}\equiv 0$ mod $p$.

To complete the proof, let $\gamma$ be an
integer $\gg 0$; we can still modify the coefficients as
${c_{jl}''=c_{jl}+ p^{\gamma}t_{j}}$. It is enough to notice that,
setting ${d={\rm lcm}\{d_i\,|\,d_{i}\not\equiv 0\;\mbox{ mod }
p\}}$, then $d$ and $p$ are coprime and the system of
congruences:
\begin{equation}
c_{jl}+ p^{\gamma}t_{j}\equiv 0\hspace{10mm}{\rm mod}\,d
\hspace{15mm}\forall\,j
\end{equation}
admits a solution. So we can assume that $c_{jl}''\equiv 0$ mod $d$,
and the proof is complete. \hfill\qed

\vspace{3mm}
To any decomposition (\ref{decteor}) as in Lemma \ref{coeffinK}, we
associate a cohomology class in $H^2(X,K)$:
\begin{Definition}\label{defxi}
Given a decomposition (\ref{decteor}) as in Lemma \ref{coeffinK},
consider the map:
\begin{equation}
\begin{array}{ccl}
{\bf Z}^{q}&\to &N\\
(x_1,\ldots x_q) &\mapsto&\sum_{l=1}^q\,x_l\underline{c}_l=
(\sum_{l=1}^q\,x_lc_{jl})_{j=1,\ldots k}
\end{array}
\end{equation}
and denote by ${\Theta:{\bf Z}^{q}\to K}$ its composition with
the projection ${N\to K}$ (cf. Prop.\ref{top}, c)). Then, set:
\begin{equation}
\xi=\Theta_*([M_1], \ldots [M_q])\,,
\end{equation}
where ${\Theta_*:H^2(X,{\bf Z}^q)\cong\bigoplus^q H^2(X,{\bf Z})\to
H^2(X,K)}$ is the map induced in cohomology by $\Theta$ and
$[M]$ is the Chern class of a divisor $M$ on $X$.
\end{Definition}

\vspace{3mm}
We briefly recall some facts about quotients by a properly
discontinuous group action (see for instance \cite{kn:Mu}, Appendix to
section 1, \cite{kn:Grot}, ch. 5).

Let $\tilde{X}$ be a simply connected variety, let $\Gamma$
be a group acting properly and discontinuously on
$\tilde{X}$ and let ${p:\tilde{X}\rightarrow X=\tilde{X}/\Gamma}$
be the projection onto the quotient. Consider the following
two functors:
$$
\begin{array}{l}
M\stackrel{F}{\longrightarrow} M^{\Gamma} ,\,{\rm for}\,M
\,{\rm a}\, \Gamma \mbox{-module}\\
{\cal F}\stackrel{H}{\longrightarrow}
H^0 (\tilde{X},p^{*}{\cal F}),\,{\rm for}\,{\cal F}
\,{\rm a\,locally\,constant\,sheaf \, on}\, X\,.
\end{array}
$$
The spectral sequence associated to the functor $F\circ H$
yields in this case the exact sequence of cohomology
group:
\begin{equation}
\label{spseq}
0\longrightarrow H^{2}(\Gamma,H^{0}(\tilde{X},p^{*}{\cal F}))
\longrightarrow  H^{2}(X,{\cal F}) \longrightarrow
H^{2}(\tilde{X},p^{*}{\cal F})^{\Gamma}
 \end{equation}
that will be used several times in the following and it
is natural with respect to the sheaf maps on $X$.

\begin{Theorem}\label{mtgen}
Let $X$, $Y$ be smooth projective varieties over {\bf C} of dimension
at least 2. Let ${f:Y \rightarrow X}$ be a totally ramified finite
abelian cover branched on a divisor with flexible and ample
components ${\{D_{j}\}_{j=1, \ldots k}}$. According to Prop.\ref{top},
b), the map $f$ induces a central extension:  \begin{equation}
\label{extgrouteogen}
0\longrightarrow K \longrightarrow\pi_{1}(Y)
\stackrel{f_{*}}{\longrightarrow} \pi_{1}(X) \longrightarrow 1
\end{equation}
Denote by ${c(f)\, \in \,H^{2}(\pi_{1}(X), K)\subseteq H^{2}(X, K)}$
the cohomology class classifying the extension (\ref{extgrouteogen}).

Let:
\begin{equation} \label{decteor}
\left(\begin{array}{c}
D_{1}\\
\vdots\\
D_{k}
\end{array}\right)\;\;\equiv\;\;C\,
\left(\begin{array}{c} M_{1}\\
\vdots\\
M_{q}
\end{array}\right)
\end{equation}
be a decomposition as in Lemma \ref{coeffinK} and let $\xi\,\in\,
H^{2}(X,K)$ be the class defined in Def.\ref{defxi}.

In this notation, one has:
\begin{equation}
c(f)=\xi;
\end{equation}
in particular, the class $\xi$ does not
depend on the chosen decomposition.
\end{Theorem}

{\sc Proof.} It is enough to show that $\xi$ and $c(f)$ admit
cohomologous representatives. This can be done in three steps.

\vspace{3 mm}
{\sc Step I}: we compute a cocycle representing $c(f)\,\in\,H^2(X,K)$.

We start by choosing suitable trivializations of the line bundles
that appear in the computation.

Set $\Gamma = \pi_{1}(X)$ and $\tilde{\Gamma} = \pi_{1}(Y)$. Let
$\{U_{r}\}$ be a sufficiently fine cover of $X$ such that $\Gamma$
acts transitively on the set of connected components of
$\pi^{-1}(U_{r})$, $\forall$ $r$. If we fix a component $V_{r}$ of
$\pi^{-1}(U_{r})$, then ${\pi^{-1}(U_{r}) = \cup_{\gamma \in
\Gamma}\,\gamma(V_{r})}$; for every ${\gamma \in \Gamma}$ we write:
\begin{equation}
\gamma(V_{r}) = V_{(\gamma, r)}
\end{equation}
and, in particular: $V_{(1, r)}=V_{r}$.

Such a covering has the following properties:
\begin{list}%
{\alph{alph})}{\usecounter{alph}}
\item For every $(r, s)$ such that ${U_{r}\cap U_{s} \neq \emptyset}$,
there exists a unique element ${\beta (r, s) \in \Gamma}$ such that:
\begin{equation}
V_{(1, r)} \cap V_{(\beta (r, s), s)} \neq \emptyset\;.
\end{equation}
\item If ${U_{r}\cap U_{s} \neq \emptyset}$, then $V_{(\gamma, r)}$
and $V_{(\gamma \beta (r, s), s)}$ have nonempty intersection.
\item Since $\pi$ is a local homeomorphism, if
${U_{r}\cap U_{s}\cap U_{t} \neq \emptyset}$, then:
\begin{equation}
\emptyset\;\neq\; V_{(\beta (r, s), s)}\cap V_{(\beta (r, t), t)}.
\end{equation}
Hence the following relation is satisfied for every
${U_{r}\cap U_{s}\cap U_{t} \neq \emptyset}$:
\begin{equation}
\beta (r, t) = \beta (r, s) \beta (s, t)\;.
\end{equation}
In particular: $\beta (s, r) = \beta (r, s)^{-1}$.
\end{list}

\noindent
For later use, we set:
\begin{equation}
V_{(\alpha, r, s)} = \alpha(V_{(1, r)}\cap V_{(\beta (r, s), s)})
= V_{(\alpha, r)}\cap V_{(\alpha\beta (r, s), s)}
\end{equation}
for every $\alpha \in \Gamma$ and for every $(r, s)$ such that
${U_{r}\cap U_{s} \neq \emptyset\,}$.

For every $r$ and for every $j = 1, \ldots k$, we choose a local
generator $w^{j}_{r}$ for ${{\cal O}_{X}(-D_{j})}$
on $U_{r}$ (we ask that $w^{j}_{r}$ is a local equation for $D_j$) and for
every pair $(r, s)$ such that ${U_{r}\cap U_{s} \neq \emptyset}$ we
write:
\begin{equation}\label{coc-Dj}
w^{j}_{r} = k^{j}_{(r, s)}\,w^{j}_{s}\;\;{\rm on}\,U_{r}\cap U_{s}\,.
\end{equation}

Now we apply to ${\tilde{f\,}:\tilde{Y}\rightarrow \tilde{X}}$ the
analysis of section 2 , most of which can be easily extended to the
case of analytic spaces. One has: \begin {equation} \label{split}
\tilde{f\,}_{*}({\cal O} _{\tilde{Y}})=
\bigoplus_{\tilde{\chi} \in \tilde{G}^{*}}L_{\tilde{\chi}}^{-1}
\end{equation}

Each element of the group $\tilde{G}$ can be interpreted as an
automorphism of the sheaf $\tilde{f\,}_{*}({\cal O} _{\tilde{Y}})$. In
particular, by duality, the elements of $K\subseteq \tilde{G}$ are
characterized by the property that they induce the identity on the
subsheaf $L_{\chi}$ for every
${\chi\,\in\,G^*\subseteq\tilde{G}^*}$.

Let $\tilde{\chi}_{1}, \ldots \tilde{\chi}_{h} \in \tilde{G}^{*}$ be
such that $\tilde{G}^{*}$ is isomorphic to the direct sum of the
cyclic subgroups generated by the $\tilde{\chi}_{i}$'s, and let
$\tilde{d_i}$ be the order of $\tilde{\chi}_{i}$, $i = 1, \ldots h$.
Let $\tilde{D_{i}}$ be the inverse image of $D_{i}$ via the universal
covering map ${\pi :\tilde{X}\to X}$, as before. If
${\tilde{\chi}_i=\sum_{j=1}^{k}\,\tilde{a}_{ij}\psi_j}$, with ${0\leq
\tilde{a}_{ij} <m_j}$, the system (\ref{eqstr}) yields in this case:
\begin{equation} \label{eqstr1}
\tilde{d_i} L_{\tilde{\chi}_{i}} \equiv
\sum_{j=1}^{k}\frac{\tilde{d_i}\tilde{a}_{ij}}{m_{j}} \;\tilde{D}
_{j}\;\;\;\;i= 1, \ldots h\,.
\end{equation}
So it is possible to choose local generators $\tilde{z}^{i}_{(\alpha,
r)}$ for $L_{\tilde{\chi}_{i}}^{-1}$ on $V_{(\alpha, r)}$ such
that for every ${\alpha \in \Gamma}$ and for every pair $(r, s)$ with
${U_{r}\cap U_{s} \neq \emptyset}$ one has:
\begin{equation}\label{relforz}
\left(\tilde{z}_{(\alpha, r)}^{i}\right)^{\tilde{d_i}}=
\prod_{j=1}^k\left(w_r^j\right)^{\frac{\tilde{d_i}\tilde{a}_{ij}}{m_j}}
\,.
\end{equation}
Writing:
 \begin{equation}
\tilde{z}_{(\alpha, r)}^{i} = \tilde{h}^{i}_{(\alpha, r,s)}\,
\tilde{z}^{i}_{(\alpha \beta (r, s), s)} \;\;\;\; {\rm
on}\,V_{(\alpha, r, s)}
\end{equation}
we have:
\begin{equation}
\label{rela}
(\tilde{h}^{i}_{(\alpha, r,
s)})^{\tilde{d_i}}\;=\;\prod_{j=1}^{k}(k^{j}_{(r, s)})
^{\frac{\tilde{d_i}\tilde{a}_{ij}}{m_{j}}}\;.
\end{equation}
and the cocycle condition for $\tilde{h}^{i}_{(\alpha, r,
s)}$, that will be often used later on, yields
the relation:
\begin{equation}
\label{coch}
1=\tilde{h}^{i}_{(\alpha, r, s)}\tilde{h}^{i}_{(\alpha\beta(r,s), s,
t)} \tilde{h}^{i}_{(\alpha\beta(r,t), t, r)}\;
\hspace{15mm}\forall\,\alpha\,\in\,\Gamma,\forall\,i,\forall\,
r,s,t\mbox{ with }U_r\cap U_s \cap U_t\neq\emptyset. \end{equation}

We observe that the generator $\tilde{z}^{i}_{(\alpha,r)}$ is
determined by (\ref{relforz}) only up to a constant of the form
$exp2\pi \sqrt{-1} \,(u^i_{(\alpha,r)}/\tilde{d_i})$ with
${u^i_{(\alpha,r)}\in {\bf Z}}$. Moreover, according to
(\ref{Lchigen}), every choice of local generators $w^j_r$ for ${\cal
O}_X(-D_j)$ and $\tilde{z}_{(\alpha, r)}^{i}$ for
$L_{\tilde{\chi}_i}^{-1}$ induces a choice of
local generators for $L_{\tilde{\chi}}^{-1}$,
${\forall\,\tilde{\chi}\,\in\,\tilde{G}^*}$, by the rule:
\begin{equation}\label{ztildechi}  \tilde{z}_{(\alpha,
r)}^{\tilde{\chi}}=\prod_{i=1}^n\left( \tilde{z}_{(\alpha, r)}^{i}
\right)^{b_{\tilde{\chi},i}}
\prod_{j=1}^k\,(w^j_r)^{-q^{\tilde{\chi}}_{j}}\,\hspace{15mm} \mbox{if
}\,\tilde{\chi}=\sum_{i=1}^h b_{\tilde{\chi},i}\tilde{\chi}_i,\;\;
0\leq b_{\tilde{\chi},i}<\tilde{d_i}\,.
 \end{equation}
where $q^{\tilde{\chi}}_j$ denotes the integral part of the real
number $\sum_{i=1}^h b_{\tilde{\chi},i} \,\tilde{a}_{ij}$.

Let now $\chi_1, \ldots \chi_n$ be a set of generators for $G^*$ such
that $G^*$ is the direct sum of the cyclic subgroups generated by the
$\chi_v$'s and the order $d_v$ of $\chi_v$ is a power of a prime
number, $v=1, \ldots n$. We recall that ${G^*\subseteq \tilde{G}^*}$
and,
${\forall\,\chi\,\in\,G^*}$, the corresponding eigensheaf $L_{\chi}$ is
a pullback from $X$. We write
${\chi_v=\sum_{i=1}^n b_{vi}\tilde{\chi}_i\,\in\,\tilde{G}^*}$
(${0\leq b_{vi}<\tilde{d_i}}$) and $q^{\chi_v}_j=q^{v}_j$; the
corresponding local generator for ${L_{\chi_v}^{-1}}$ chosen in
(\ref{ztildechi}) is:
\begin{equation}\label{zG*}
z_{(\alpha, r)}^{v}=\prod_{i=1}^h\left(
\tilde{z}_{(\alpha, r)}^{i}
\right)^{b_{vi}}\prod_{j=1}^k\,(w^j_r)^{-q^{v}_j} \;;
\end{equation}
we show that, for a suitable choice of the $\tilde{z}_{(\alpha,
r)}^{i}$, we can assume that the expression in (\ref{zG*}) is
independent from $\alpha$. In fact, using the characteristic equations
of the cover $f$, one can choose a local base $y^v_r$ of
$L_{\chi_v}^{-1}$ on $V_{(\alpha, r)}$ that does not depend on
$\alpha$ and satisfies the relation:
\begin{equation}\label{eqstrsuX}
\left(y^v_r\right)^{d_v}=\prod_{j=1}^k\left(w^j_r
\right)^{\frac{a_{vj} d_v}{m_j}}\,.
\end{equation}

Since $\sum_{i=1}^h b_{vi}\tilde{a}_{ij}=q^v_jm_j+a_{vj}$
$\forall\,j=1,\ldots k$, the two local generators $y^v_r$ and
$z_{(\alpha, r)}^{v}$ on $V_{(\alpha,r)}$ differ by a $d_v$-th
root of unity, that we denote by $exp( 2 \pi
\sqrt{-1}\,(x^v_{(\alpha,r)}/d_v))$. If we multiply
$\tilde{z}_{(\alpha, r)}^{i}$ by  $exp( 2 \pi
\sqrt{-1}(u^i_{(\alpha,r)}/\tilde{d_i}))$, then
${(x^v_{(\alpha,r)}/d_v)}$ becomes
${(x^v_{(\alpha,r)}/d_v)+\sum_{i=1}^h
(b_{vi}/\tilde{d_i})\,u^i_{(\alpha, r)}}$. Hence, we only need to
solve the linear system of congruences, $\forall\,(\alpha, r)$:
\begin{equation} \sum_{i=1}^h\frac{d_v
b_{vi}}{\tilde{d_i}}u^i_{(\alpha, r)}\equiv  x^v_{(\alpha,
r)}\;\;\;\mbox{mod}\,d_v\hspace{20 mm}v=1, \ldots n; \end{equation}
since we assume that $d_v$ is a power of a prime number, this system
admits a solution according to the Chinese Remainder's Theorem and
Lemma \ref{Lemma}.

So we can assume that the expression $z_{(\alpha, r)}^{v}$ in
(\ref{zG*}) does not depend on $\alpha$ and it is the pullback of a
local generator of the corresponding eigensheaf on $X$: we write
$z_{r}^{v}=z_{(\alpha, r)}^{v}$. For later use, we define the
corresponding cocycle $h^{v}_{(r,s)}$ ($v=1,\ldots n$) by the rule:
\begin{equation}\label{cochv}
z_{r}^{v}=h^{v}_{(r,s)}z_{s}^{v} \hspace{15 mm}\mbox{on }U_r\cap
U_s
\end{equation}
and we observe that, according to (\ref{eqstrsuX}), the following
relation holds:
\begin{equation}\label{relhv}
(h^{v}_{(r,s)})^{d_v}=\prod_{j=1}^k(k^j_{(r,s)})^{\frac{a_{vj}
d_v}{m_j}}\,\hspace{15mm}\mbox{if
}\chi_v=\sum_{j=1}^k a_{ij}\psi_j, \;\;0\leq a_{ij}<m_j.
\end{equation}

\vspace{3 mm}
In order to compute the class of the extension
(\ref{extgrouteogen}), for every ${\gamma \in \Gamma}$ we choose a
lifting ${\tilde{\gamma} \in \tilde{\Gamma}}$.
By Lemma \ref{commut}, the induced map ${\tilde{\gamma}_{*}:
\tilde{f\,}_{*} {\cal O}_{\tilde {Y}} \to \tilde{f\,}_{*}{\cal
O}_{\tilde {Y}}}$ is a  ${\cal O}_{\tilde {X}}$-algebra  isomorphism
lifting ${\gamma : \tilde{X} \to \tilde{X}}$; in terms of the chosen
trivializations we may write:
\begin{equation}
\tilde{z}^{i}_{(\alpha ,r)} \stackrel{\tilde{\gamma}_{*}}{\mapsto}
\sigma^{i,\gamma}_{(\alpha,r)}\tilde{z}^{i}_{(\gamma\alpha
,r)}\;\;\;\;\;\; \;\forall\gamma , \alpha \in \Gamma , i=1,\ldots h
\end{equation}
for a suitable choice of a $\tilde{d}_{i}$-th root of unity
$\sigma^{i,\gamma}_{(\alpha,r)}$, $i=1, \ldots h$.

For later use, we write down the transition relations
for the constants $\sigma^{i,\gamma}_{(\alpha,r)}$. Let $s,t$ be such that
${U_{s}\cap U_{t} \neq \emptyset}$ and let $\alpha,\gamma \in\Gamma$;
then, for $i=1,\ldots h$:
\begin{equation}
\label{rel}
\sigma^{i,\gamma}_{(\alpha ,s)}\tilde{h}^{i}_{(\gamma\alpha ,s,t)}=
\sigma^{i,\gamma}_{(\alpha\beta (s,t),t)}\left(\tilde{h}^{i}_{(\alpha
,s,t)}\circ \gamma^{-1}\right)\;\;\;\;\;\;\;{\rm
on}\;V_{(\gamma\alpha, s,t)}.
\end{equation}

We now exploit the action of the chosen elements in $\tilde{\Gamma}$ on
$\tilde{f\,}_{*} {\cal O}_{\tilde {Y}}$ in order to compute the class
${c(f)\in H^{2}(\Gamma ,K)}$
associated to the extension (\ref{extgrouteogen}), and its image
${c(f)\in H^{2}(X,K)}$.

For any given
$\delta,\gamma\in\tilde{\Gamma}$, the action of
$(\widetilde{\delta\gamma})_{*}^{-1}\tilde{\delta}_{*}
\tilde{\gamma}_{*}$ on $L^{-1}_{\tilde{\chi}_{i}}$, $i = 1,
\ldots h$, is described with respect to the chosen trivializations by:
\begin{equation} \label{comm}
\tilde{z}^{i}_{(\alpha ,r)} \mapsto
\left(\sigma^{i,(\delta\gamma)}_{(\alpha,r)}\right)^{-1}
\sigma^{i,\delta}_{(\gamma\alpha ,r)} \sigma^{i,\gamma}_{(\alpha ,r)}
\tilde{z}^{i}_{(\alpha ,r)}\;\;\;\;\;\;\;\;\;\;\forall
r,\;\forall\alpha\in\Gamma, \; i=1,\ldots h.
\end{equation}

Since (\ref{comm}) represents a line bundle automorphism given by a
root of the unity, the expression does not depend on $(\alpha ,r)$ by
the connectedness of $\tilde{X}$: therefore, we may set $\alpha = 1$.
So, the class ${c(f)\in H^{2}(\Gamma
,K)}$ is represented by the cocycle:
\begin{equation}
c(f)(\delta ,\gamma )=\left(
\left(\sigma^{i,(\delta\gamma)}_{(1,r)}\right)^{-1}
\sigma^{i,\delta}_{(\gamma,r)}
\sigma^{i,\gamma}_{(1,r)}\right)_{i=1,\ldots h}\;
\;\;\;\;\;\forall r,
\end{equation}
where an element of $K\subseteq \tilde{G}$ is represented by
its coordinates with respect to the
basis dual to $\{\tilde{\chi}_{1}, \ldots \tilde{\chi}_{h}\}$.

According to (\cite{kn:Mu}, page 23), the class
${c(f)\in H^{2}(X,K)}$ is represented on $V_{(1,r,s)}\cap V_{(1,
r,t)}$ by the cocycle:
\begin{equation}
\label{uffa}
c(f)_{r,s,t}=c(f)(\beta (r,s),\beta (s,t))=
\left(
\left(\sigma^{i,\beta (r,t)}_{(1,p)}\right)^{-1}
\sigma^{i,\beta (r,s)}_{(\beta (s,t),p)}
\sigma^{i,\beta (s,t)}_{(1 ,p)}
\right)_{i=1,\ldots h}\;
\;\;\;\;\;\;\forall p
\end{equation}
for $r,s,t$ such that
${U_{r}\cap U_{s}\cap U_{t}\neq\emptyset}$.

We set $p=t$ and, by the relation (\ref{rel}), we rewrite (\ref{uffa})
as follows:
\begin{equation}\label{coccf2}
c(f)_{r,s,t}=\left(
\left(\sigma^{i,\beta (r,t)}_{(1,p)}\right)^{-1}
\sigma^{i,\beta (t,r)}_{(1,r)}
\sigma^{i,\beta (s,t)}_{(1 ,t)}
 \tilde{h}^{i}_{(\beta (r,s),s,t)}
\left( \tilde{h}^{i}_{(1,s,t)}\circ \beta(s,r)\right)^{-1}
\right)_{i=1, \ldots h}\,;
\end{equation}
this shows that $c(f)_{r,s,t}$ differs from the following cocycle (that
we still denote by $c(f)_{r,s,t}$ by abuse of notation):
\begin{equation}\label{coccfnuovo}
c(f)_{r,s,t}=\left(
\tilde{h}^{i}_{(\beta (r,s),s,t)}
\left( \tilde{h}^{i}_{(1,s,t)}\circ \beta(s,r)\right)^{-1}\right)_{i=1, \ldots
h}\, \end{equation}
by the coboundary of the cochain:
\begin{equation}
g_{r,t}=\left(\sigma^{i,\beta(r,t)}_{(1,t)}\right)_{i=1, \ldots h}\,.
\end{equation}
The cochain $g_{r,t}$ actually takes values in $K$: in fact, it is
enough to check $g_{r,t}$ acts trivially on the eigensheaves
corresponding to the chosen generators $\chi_v$ of $G^*$. This
follows easily by the prevous choices since the action of $\beta(r,t)$
on $L_{\chi_v}^{-1}$ is given locally by
${\prod_{i=1}^h\left(\sigma^{i,\beta(r,t)}_{(1,t)}\right)^{b_{vi}}}$.

\vspace{3mm}
{\sc Step II:} we compute a cocycle representing
$\xi\,\in\,H^2(X,K)$.

For every $r$ and for every $l = 1, \ldots q$, we choose a local
generator $y^{l}_{r}$ for ${\cal O}_X(-M_{l})$ on $U_{r}$; if $M_{l}$
has finite order $e$, then we require:
\begin{equation}
\label{tors}
\left(y^{l}_{r}\right)^{e}=1\,.
\end{equation}

We set $m={\rm lcm}\,\{m_j\}_{j=1, \ldots k}$. For every pair of
indices $(r, s)$ such that ${U_{r}\cap U_{s} \neq \emptyset}$ we write:
\begin{equation}
y^{l}_{r} = \mu^{l}_{(r, s)}\,y^{l}_{s}\;\;{\rm on}\,U_{r}\cap U_{s}
\end{equation}
and we choose a $m$-th root $\hat{\mu}^{l}_{(r, s)}$ of
$\mu^{l}_{(r, s)}$ in such a way that
${\hat{\mu}^{l}_{(s, r)}=(\hat{\mu}^{l}_{(r,
s)})^{-1}}$. Then, as in \cite{kn:Cato}, (2.45), one sees that
the image of the class of $-M_{l}$ in $H^2(X,{\bf Z}/m_l{\bf Z})$ is
represented on  ${U_{r}\cap U_{s}\cap U_{t}}$ by the cocycle
${\left(\hat{\mu}^{l}_{(r, s)}
\hat{\mu}^{l}_{(s, t)}
(\hat{\mu}^{l}_{(r, t)})^{-1} \right)^{m/m_l}}$, $l=1, \ldots q$. We
conclude that the class ${\xi=\Theta_*([M_1], \ldots [M_q])}$ is
represented on ${U_{r}\cap U_{s}\cap U_{t}}$ by:
\begin{equation}\label{cocc}
\xi_{r,s,t}=\left(\prod_{j=1}^k\prod_{l=1}^{q}\,\left(\hat{\mu}^{l}_{(r,
s)} \hat{\mu}^{l}_{(s, t)}
\hat{\mu}^{l}_{(t, r)}\right)^{-(m/m_j)
c_{jl}\tilde{a}_{ij}}\right)_{i=1, \ldots h}
\end{equation}

\vspace{3mm}
{\sc Step III:} we show that ${\xi=c(f)}$.

We remark that, according to (\ref{decteor}), the cocycle $k^{j}_{(r,
s)}$ in (\ref{coc-Dj}) representing ${{\cal O}_{X}(-D_{j})}$ ($j=1,
\ldots k$) and the cocycles $\mu^{l}_{(r,s)}$ representing $-M_l$
($l=1, \ldots q$) are related as follows:
\begin{equation}
k^{j}_{(r, s)}=\prod_{l=1}^q\left(\mu^{l}_{(r,
s)}\right)^{c_{jl}}\frac{f^j_r}{f^j_s}
\end{equation}
for suitable nowhere vanishing holomorphic functions $f_r$ on
$U_{ r}$. For every $j=1, \ldots k$ and every $r$, we choose a $m$-th
root $\hat{f}^j_r$ of $f^j_r$ on $U_{ r}$; then, the
expression:
\begin{equation}\label{tildekmu}
\hat{k}^j_{(r,s)}=\prod_{l=1}^{q}\,\left(\hat{\mu}^{l}_{(r,
s)}\right)^{c_{jl}(m/m_j)}
\left(\frac{\hat{f\,}^j_r}{\hat{f\,}^j_s}\right)^{(m/m_j)}
\end{equation}
is a $m_j$-th root of the
cocycle $k^j_{(r,s)}$ and, as before, the product
$\prod_{l=1}^{q}\,\left(\hat{\mu}^{l}_{(r, s)} \hat{\mu}^{l}_{(s,
t)} \hat{\mu}^{l}_{(t, r)}\right)^{ c_{jl}(m/m_j)}$ yields a cocycle
representing the image of the class of $-D_j$ in $H^2(X,{\bf
Z}/m_j{\bf Z})$. In this notation, by (\ref{tors}), we rewrite
as follows the cocycle in (\ref{cocc}) representing the class $\xi$:
\begin{equation}\label{cocc1}
\xi_{r,s,t}=\left(\prod_{j=1}^k
\left( \hat{k}^j_{(r,s)}\hat{k}^j_{(s,t)}\hat{k}^j_{(t,r)}
\right)^{-\tilde{a}_{ij}}\right)_{i=1, \ldots, h}\,.
\end{equation}

Let $\epsilon = exp(\frac{2\pi \sqrt{-1}}{m})$. Then, by
(\ref{rela}), one has:
\begin{equation}\label{defq}
\tilde{h}^{i}_{(\alpha, r, s)}\;=\;\prod_{j=1}^{k}(\hat{k}^{j}_{(r,s)})
^{\tilde{a}_{ij}}\;
\epsilon^{-q^{i}_{(\alpha,r,s)}}
\end{equation}
where $q^{i}_{(\alpha,r,s)}$ is an integer, multiple of
$m/\tilde{d_i}$, and (\ref{coccfnuovo}) may be rewritten as:
\begin{equation}
c(f)_{r,s,t} =
(\epsilon^{q^{i}_{(1,s,t)}-q^{i}_{(\beta(r,s),s,t)}})
_{i=1,\ldots h}\;\;.
\end{equation}
{}From the cocycle condition (\ref{coch}) for
$\tilde{h}^{i}_{(\alpha , r, s)}$, it follows:
\begin{equation}
\xi_{r,s,t}=(\epsilon^{-q^{i}_{(\alpha
,r,s)}-q^{i}_{(\alpha\beta (r,s),s,t)}
-q^{i}_{(\alpha\beta(r,t),t,r)}})_{i=1,\ldots h}\;\;\;\forall\, r,s,t,
\forall\,\alpha\in\Gamma\,.
\end{equation}

In particular, for $\alpha = 1$, one gets:
\begin{equation}
\xi_{r,s,t}=(\epsilon^{-q^{i}_{(1,r,s)}-q^{i}_{(\beta
(r,s),s,t)} -q^{i}_{(\beta (r,t),t,r)}})_{i=1,\ldots h}\;\;.
\end{equation}
So, one has:
\begin{equation}
c(f)_{r,s,t} = \xi_{r,s,t}
(\epsilon^{q^{i}_{(1,r,s)}+q^{i}_{(1,s,t)}
+q^{i}_{(\beta(r,t),t,r)}})_{i=1,\ldots h}\;\;.
\end{equation}
By the definition of $q^{i}_{(\alpha,r,s)}$, this equality can
be rewritten as follows:
\begin{equation}\label{cobord}
c(f)_{r,s,t} = \xi_{r,s,t}
(\epsilon^{q^{i}_{(1,r,s)}-q^{i}_{(1,r,t)}
+q^{i}_{(1,s,t)}})_{i=1,\ldots h}\,.
\end{equation}

To complete the proof of the theorem, we show that we can
choose the $m$-th root $\hat{f\,}_j$  of $f_j$ ($j=1, \ldots k$) so
that:
\begin{equation}
\underline{q}_{(1,r,s)}=\left(\epsilon^{q^{i}_{(1,r,s)}}
\right)_{i=1,\ldots h} \mbox{ is an element of } K, \;\forall\,(r,s).
\end{equation}
Let $\chi_v$ one of the chosen generators of $G^*$. According to
(\ref{zG*}), (\ref{relhv}) and (\ref{defq}), the action of
$\underline{q}_{(1,r,s)}$ on $L_{\chi_v}^{-1}$ is given by a $d_v$-th
root of unity, that we denote by $exp({2\pi
\sqrt{-1}\,\frac{x_v}{d_v}})$ (for a suitable integer $x_v$). We want
to show that we can assume that $x_v\equiv 0$ mod $d_v$, $v=1, \ldots
n$.

We observe that:
\begin{equation}\label{cobinK}
exp({2\pi \sqrt{-1}\,\frac{x_v}{d_v}})=
\prod_{i=1}^h \epsilon^{q^{i}_{(1,r,s)}b_{vi}}=
(h^{v}_{(1,r,s)})^{-1}\prod_{j=1}^k(\hat{k}_{(r,s)}^j)^{a_{vj}}
\end{equation}
and we compute the right-hand side of
(\ref{cobinK}).  By (\ref{tildekmu}), one must have:
\begin{equation}
\begin{array}{cl}
h^{v}_{(1,r,s)}&=exp({-2\pi
\sqrt{-1}\,\frac{x_v}{d_v}})\prod_{j=1}^k
(\hat{k}^j_{(r,s)})^{a_{vj}}\\
&= exp({-2\pi
\sqrt{-1}\,\frac{x_v}{d_v}}) \prod_{l=1}^q
(\hat{\mu}^l_{(r,s)})^{m\sum_{j=1}^k \frac{c_{jl}a_{vj}}{m_j}}
\prod_{j=1}^k\left(\frac{\hat{f\,}^j_r}{\hat{f\,}^j_s}
\right)^{a_{vj}(m/m_j)}\\
&= exp({-2\pi
\sqrt{-1}\,\frac{x_v}{d_v}}) \prod_{l=1}^q
(\hat{\mu}^l_{(r,s)})^{\frac{m}{d_v}\sum_{j=1}^k
\frac{c_{jl}d_v a_{vj}}{m_j}}
\prod_{j=1}^k\left(\frac{\hat{f\,}^j_r}{\hat{f\,}^j_s}
\right)^{a_{vj}(m/m_j)}. \end{array}
\end{equation}
On the other hand, as in Lemma
\ref{coeffinK}, we write $L_{\chi_v}\equiv \sum_{l=1}^q
\lambda_{vl}M_l$ and we get the following relation form of cocycles on
$V_{(1,r,s)}$:    \begin{equation}
h^{v}_{(1,r,s)}=\prod_{l=1}^q
(\mu^l_{(r,s)})^{\lambda_{vl}}
\frac{\varphi^v_{r}}{\varphi^v_{s}}
\end{equation}
for suitable nowhere vanishing holomorphic functions $\varphi^v_{r}$
on $U_{r}$. According to Lemma \ref{Lemma} and to (\ref{tors}), we can
then assume that in the previous equations one has:
\begin{equation}
 \prod_{l=1}^q
(\hat{\mu}^l_{(r,s)})^{\left(\frac{m}{d_v}\sum_{j=1}^k
\frac{c_{jl}d_v a_{vj}}{m_j}\right)-d_v \lambda_{vl}}=1
\end{equation}
so that one gets:
\begin{equation}
\label{xv}
exp({2\pi \sqrt{-1}\,\frac{x_v}{d_v}})=
\prod_{j=1}^k\left(\frac{\hat{f\,}^j_r}{\hat{f\,}^j_s}
\right)^{a_{vj}(m/m_j)}\frac{\varphi^v_s}{\varphi^v_r}\,.
\end{equation}
We observe that we may assume that:
\begin{equation}
\varphi^v_r= exp(2\pi \sqrt{-1}\,\frac{t^v_r}{d_v})
\prod_{j=1}^k(\hat{f\,}^j_r)^{a_{vj}(m/m_j)}
\end{equation}
for suitable integers $t^v_r$; hence the equation
(\ref{xv}) gives:
\begin{equation}
\frac{x_v}{d_v}- \frac{t^v_s}{d_v}+\frac{t^v_r}{d_v}\;\in\;{\bf Z}.
\end{equation}
If we replace $\hat{f\,}^j_r$ by $exp(2\pi
\sqrt{-1}\,\frac{s^j_r}{m})\hat{f\,}^j_r$, then
$\frac{t^u_r}{d_v}$ is replaced by
$\frac{t^u_r}{d_v}+\sum_{j=1}^k\frac{a_{vj}s_j}{m_j}$. Therefore, we
need to solve the system:
\begin{equation}
\sum_{j=1}^k\frac{a_{vj}s^j_r}{m_j}\equiv t^v_r \hspace{10mm}
\mbox{mod}\,d_v \hspace{15mm}v=1,\ldots n.
\end{equation}
Since this is possible according to Lemma \ref{Lemma} and the Chinese
Remainder's Theorem, the proof is complete.
\hfill\qed

\begin{Remark}\label{dipdaLchi}
The cohomology class of the extension (\ref{extgrouteogen}) of the
fundamental groups depends on the choice of the solution
$\{L_\chi\}$ of the characteristic relations (\ref{eqstr}) for the
covering $f$. Moreover, covers corresponding to different
solutions $\{L_\chi\}$ may not be homeomorphic.

{\rm This is shown, for instance, by the following class
of examples. Denote by $e_i$ the standard generators of the group
$({\bf Z}/4{\bf Z})^3$ and let $G$ be the quotient of $({\bf Z}/4{\bf
Z})^3$ by the subgroup generated by ${2e_1+2e_2+2e_3}$. Let now $X$
be a smooth projective surface such that ${H^2(\pi_1(X), {\bf Z}/2
{\bf Z})\neq 0}$ and Pic($X$) has a 2-torsion element $\eta$ whose
class in ${H^2(X, {\bf Z}/2{\bf Z})}$ is non zero. Fix a very ample
divisor $H$ on $X$ and choose suitable divisors $D_i$ ($i=1,2,3$) such
that $D_i\equiv 4H$ and the $D_i$'s are in general position. Then
there exists a smooth abelian $G$-cover ${f:Y\to X}$ ramified on the
$D_i$'s ($i=1,2,3$), with inertia subgroup $G_i=<e_i>={\bf Z}/4{\bf
Z}$ and character $\psi_i$ dual to $e_i$, respectively. In fact,
taking the characters $\chi_1=\psi_1+3\psi_3$,
$\chi_2=\psi_2+3\psi_3$, $\chi_3=2\psi_3$ as generators of $G^*$, the
characteristic relations (\ref{eqstr}) of the cover $f$ are:
\begin{equation} \label{eqstrex}
\left\{ \begin{array}{lll}
4L_1&\equiv&D_1+3D_3\\ 4L_2&\equiv&D_2+3D_3\\ 2L_3&\equiv&D_3
\end{array}
\right.
\end{equation}
and admit, in particular, the solution $L_1= L_2=
L_3= 2H$. Under these hypotheses, $L_3$ generates the
subgroup $<D_i,L_\chi>$ ($i=1,2,3,\chi\,\in\,G^*$) of Pic($X$) and the
decomposition $D_i\equiv 2L_3$ has the properties requested in
Prop.\ref{coeffinK}. According to Prop.\ref{top} and Thm.\ref{mtgen},
since the pull back $\tilde{D}_i$ of $D_i$ under the universal cover
$\tilde{X}$ of $X$ is 2-divisible, $\forall\,i$, then the map $f$
induces a central extension of the form:
\begin{equation}
0\to {\bf Z}/2{\bf Z}\to \pi_1(Y)\to \pi_1(X)\to 1
\end{equation}
and the cohomology class of
this extension in $H^2(\pi_1(X), {\bf Z}/2{\bf Z})\subseteq H^2(X,{\bf
Z}/2{\bf Z})$ is the image $\Psi_*([L_3])$ of the Chern class of
$L_3$ under the map induced in cohomology by the standard projection
$\Psi:{\bf Z}\to {\bf Z}/2{\bf Z}$: so, this class is trivial.

Let now $\overline{Y}$ be the $G$-cover of $X$ corresponding to the
solution $\overline{L}_i=2H+\eta$ ($i=1,2,3$) of (\ref{eqstrex});
in this case the cohomology class describing $\pi_1(\overline{Y})$
is given by $\Psi_*([L_3+\eta])$ and, by the hypotheses made, it is
not trivial.

In particular, when $X$ is a
projective variety with $\pi_1(X)={\bf Z}/2{\bf
Z}$, the previous construction yields two non
homeomorphic $G$-covers $Y$, $\overline{Y}$ of $X$, branched on the
same divisor, with the same inertia subgroups and characters, such
that:
\begin{equation} \pi_1(Y)=({\bf Z}/2{\bf
Z})^2\hspace{20mm}\pi_1(\overline{Y})={\bf Z}/4{\bf
Z}\,.
\end{equation}
 }
\end{Remark}
\hfill\qed

The following theorem is an attempt to determine to what extent the
class $c(f)$ depends on the choice of the $L_\chi$'s, once the branch
divisor and the covering structure are fixed.

\begin{Theorem} \label{mt}
Same hypotheses and notation as in the statement of Thm.\ref{mtgen}.
Consider the class $c(f)\,\in\,H^{2}(\pi_1(X), K)$ associated to the
central extension (\ref{extgrouteogen}) given by the fundamental
groups and denote by  ${i(c(f)) \in
H^{2}(\pi_1(X),\tilde{G})\subseteq H^{2}(X,\tilde{G})}$ its image via
the map induced in cohomology by the inclusion (\ref{extgrouO})
$K\subseteq \tilde{G}$.

Denote by $\Phi$ the group homomorphism defined as follows:
\begin{equation}
\begin{array}{rccl}
\Phi :& {\bf Z}^{k}& \rightarrow &\tilde{G}\\
&(x_{1}, \ldots x_{k})&\rightarrow&
g_{1}^{x_{1}}\cdots g_{k}^{x_{k}}
\end{array}
\end{equation}
and by ${\Phi_{*}: H^{2}(X,{\bf
Z}^{k})  \rightarrow H^{2}(X,\tilde{G})}$ the map induced by $\Phi$
in cohomology.

Then:
\begin{equation}
\label{coc}
i(c(f)) = \Phi_{*}([D_{1}], \ldots [D_{k}])
\end{equation}
where $[\Delta]$ denotes the class of a divisor $\Delta$
on $X$ in ${H^{2}(X,{\bf Z})}$.
\end{Theorem}

\begin{Corollary}\label{mtcor}
Same hypotheses and notation as in Thm.\ref{mtgen}. Assume moreover
that the natural morphism $Hom(\pi_1(X),\tilde{G})\to
Hom(\pi_1(X),G)$, induced by the surjection ${\tilde{G}\to G}$, is
surjective. Then the map ${i:H^2(\pi_1(X),K)\to
H^2(\pi_1(X),\tilde{G})}$ is injective and the class
$\Phi_{*}([D_{1}], \ldots [D_{k}])$ in (\ref{coc}) determines
uniquely the class $c(f)\,\in\,H^2(\pi_1(X),K)$ of the extension
(\ref{extgrouteogen}) of the fundamental groups.

This happens, in particular, if $\pi_1(X)$ is torsion
free or $Hom(\pi_1(X),G)=0$ (e.g., if $\pi_1(X)$ is finite with order
coprime to the order of $G$), or the sequence ${0\to K\to
\tilde{G}\to G\to 0}$ splits.
 \end{Corollary}

{\sc Proof of Thm.}\ref{mt}. We keep the notation and the results
in Step I of the proof of Thm.\ref{mtgen}, noticing that the cocycle
$c(f)_{r,s,t}$ in (\ref{coccfnuovo}) also represents the class
$i(c(f))$ in $H^{2}(X,\tilde{G})$.

We want to write down a cocycle representing the class
${\Phi_{*}([D_{1}], \ldots [D_{k}])\in H^{2}(X,\tilde{G})}$
and to show that it represents the same cohomology class as the
cocycle in (\ref{coccfnuovo}).

We consider as before the cocycle $k^{j}_{(r,s)}$ representing
${\cal O}_{X}(-D_{j})$ in the choosen covering ${U_{r}}$ of $X$.
For every
pair of indices $r$, $s$  with ${U_{r}\cap U_{s}\neq\emptyset}$
and for every $j = 1, \ldots k$, we choose a $m_{j}$-th root
$\hat{k}^{j}_{(r,s)}$ of $k^{j}_{(r,s)}$ on ${U_{r}\cap U_{s}}$
in such a way that
${\hat{k}^{j}_{(s,r)}= (\hat{k}^{j}_{(r,s)})^{-1}}$. As before, the
image of the class of $-D_{j}$ in $H^{2}(X,{\bf Z}/m_{j}{\bf Z})$ is
represented on  ${U_{r}\cap U_{s}\cap U_{t}}$ by the cocycle
${\hat{k}^{j}_{(r,s)}\hat{k}^{j}_{(s,t)}\hat{k}^{j}_{(t,r)}}$,
${j = 1,} \ldots k$.
Then the class  ${-\Phi_{*}([D_{1}], \ldots [D_{k}])\in
H^{2}(X,\tilde{G})}$ is represented on ${U_{r}\cap U_{s}\cap U_{t}}$
by:
\begin{equation} b_{r,s,t} = (\prod^{k}_{j=1}
(\hat{k}^{j}_{(r,s)}\hat{k}^{j}_{(s,t)}\hat{k}^{j}_{(t,r)})
^{\tilde{a}_{ij}})_{i=1,\ldots h}\;\;.
\end{equation}
and we have shown in the equality (\ref{cobord}) in the proof of
Thm.\ref{mtgen}, Step III, that this cocycle represents the same
class then $c(f)_{r,s,t}$ in $H^2(X,\tilde{G})$. \hfill\qed

\begin{Remark}
\label{rem}
{\rm From Thm.\ref{mt} it follows in particular that the class
$i(c(f))\in H^{2}(\Gamma,\tilde{G})$ depends only on the class of the
$D_{j}$'s in $H^{2}(X,{\bf Z}/m_{j}{\bf Z})$ ($j = 1, \ldots k$), once
$G$ and the $g_{j}$'s are fixed. In particular, if $D_j$ is
$m_j$-divisible on $X$ ($\forall\; j=1,\ldots k$), then $i(c(f))=0$.}
\end{Remark}

\noindent
Rita Pardini\\
Dipartimento di Matematica, Universit\`a di Pisa\\
Via Buonarroti 2\\
56127 Pisa, Italy

\vspace {3mm}
\noindent
Francesca Tovena \\
Dipartimento di Matematica, Universit\`{a} di Bologna\\
Piazza di Porta S.Donato 5\\
40127 Bologna, Italy

\end{document}